\documentclass[11pt]{article}
 \usepackage[top=1in,bottom=1in,left=1.5in,right=1.5in]{geometry}
  \usepackage{amsmath,amssymb}
  \usepackage{latexsym}% defines $\Box$ for LaTeX2e
 \usepackage[dvips]{pstcol} % PSTricks
  \usepackage{pst-node}
  \usepackage[dvips]{graphicx}
  \newcommand{\C}{\mathbb{C}}
  
  \newcommand{\R}{\mathbb{R}}

  \renewcommand{\u}{\mathbf{u}}
  \renewcommand{\v}{\mathbf{v}}

  \newcommand{\x}{\mathbf{x}}
  \newcommand{\y}{\mathbf{y}}
  
  \newcommand{\0}{\mathbf{0}}

  \newcommand{\cR}{\mathcal{R}}

  \def\diag{\mathop{{\rm diag}}\nolimits}
  \newcommand{\hs}{\hspace*{\parindent}}

  \newcommand{\qed}{\hspace*{\fill} $\Box$\\}

  \newcommand{\rank}{\mathrm{rank\;}}

  \newtheorem{theo}{\bfseries \hs Theorem}[section]

  \numberwithin{equation}{section} % Automatically number equations within sections

\begin{document}

 \title{Generalized rank-constrained matrix approximations}

 \author{Shmuel Friedland
 \\Department of Mathematics, \\ Statistics and Computer
 Science\\ University of Illinois at Chicago
 \\ Chicago, Illinois 60607-7045\\ USA\\
 friedlan@uic.edu
 \and Anatoli Torokhti\\
 School of Mathematics \\and Statistics\\
 University of South Australia\\
 Mawson Lakes  SA  5095\\
 Australia \\
 Anatoli.Torokhti@unisa.edu.au }

  \date{}

 \maketitle

 \begin{abstract}
 In this paper we give an explicit solution to the rank constrained matrix approximation
 in Frobenius norm,
 which is a generalization of the classical approximation of an $m\times n$ matrix $A$
 by a matrix of rank $k$ at most.
     \\[\baselineskip] 2000 Mathematics Subject
     Classification: 15A18.
 \end{abstract}

 \section{Introduction}

 Let $\C^{m\times n}$ be set of $m\times n$ complex valued matrices, and denote by
 $\cR(m,n,k)\subseteq \C^{m\times n}$ the variety of all $m\times n$ matrices of rank
 $k$ at most.  Fix $A=[a_{ij}]_{i,j=1}^{m,n}\in \C^{m\times n}$.  Then $A^*\in\C^{n\times m}$ is
 the conjugate transpose of $A$, and
 $||A||_F:=\sqrt{\sum_{i,j=1}^{m,n}
 |a_{ij}|^2}$ is the Frobenius norm of $A$.  Recall that the \emph{singular value
 decomposition} of $A$,
 abbreviated here as \emph{SVD}, is given by  $A=U_A\Sigma_A V_A^*$, where
 $U_A\in \C^{m\times m}, V_A\in \C^{n\times n}$
 are unitary matrices,
 $\Sigma_A:=\diag(\sigma_1(A),\ldots,\sigma_{\min(m,n)}(A))\in\C^{m\times n}$
 is a generalized diagonal matrix,
  with the singular values $\sigma_1(A)\ge \sigma_2(A)\ge\ldots\ge 0$
 on the main diagonal.  The number of positive singular values of $A$ is $r$,
 which is equal to the
 rank of $A$, denoted by $\rank A$.  Let $U_A=[\u_1 \;\u_2\;\ldots\u_m],
 V_A=[\v_1 \;\v_2\;\ldots\v_n]$
 be the representations of $U,V$ in terms of their $m,n$ columns respectively.
 Then $\u_i$ and $\v_i$ are  called the \emph{left} and the \emph{right}
 singular vectors of
 $A$, respectively, that  correspond to the singular value $\sigma_i(A)$.
 Let
 \begin{equation}\label{defpalf}
 P_{A,L}:=\sum_{i=1}^{\rank A} \u_i\u_i^* \in \C^{m\times m}, \quad
 P_{A,R}:=\sum_{i=1}^{\rank A} \v_i\v_i^* \in \C^{n\times n},
 \end{equation}
 be the orthogonal projections on the range of $A$ and $A^*$ respectively.
 Denote by
 $$
 A_k:=\sum_{i=1}^k \sigma_i(A)\u_i \v_i^*\in \C^{m\times n}
 $$
  for $k=1,\ldots,\rank A$.
 For $k>\rank A$ we define $A_k:=A\;(=A_{\rank A})$. For $1\le k<\rank A$, the matrix $A_k$ is uniquely defined
 if and only if $\sigma_k(A)>\sigma_{k+1}(A)$.

 The enormous application of SVD decomposition of $A$ in pure and applied mathematics, is
 derived from the following approximation property:

 \begin{equation}\label{svdap}
 \min_{X\in \cR(m,n,k)} ||A - X||_F= ||A - A_k||_F, \quad k=1,\ldots
 \end{equation}

 The latter is known as the Eckart-Young theorem \cite{eck1}. We note
 that the work \cite{eck1} implied a number of extensions. We cite
 \cite{fri1,GV,kol1,luu1} as some recent references.
  Another application of SVD is a formula for the Moore-Penrose
 inverse $A^{\dagger}:=V_A\Sigma_A^{\dagger}U_A^*\in \C^{n\times m}$ of $A$,
 where

 \noindent
 $\Sigma_A^{\dagger}:=\diag(\frac{1}{\sigma_1(A)},\ldots,\frac{1}{\sigma_{\rank A}(A)},0,\ldots,0)
 \in \C^{n\times m}$.  See for example \cite{ben1}.

 \section{Main Result}

 Below, we provide generalizations of the classical minimal
 problem given in (\ref{svdap}).

 \begin{theo}\label{comrnkapt}  Let matrices $A\in \C^{m\times n},$ $B\in \C^{m\times p}$ and
 $C\in \C^{q\times n}$ be given.
 Then
 \begin{equation}
 \label{sol1}
 X=B^{\dagger}(P_{B,L}AP_{C,R})_k C^{\dagger}
 \end{equation}
  is a solution to the minimal problem

 \begin{equation}\label{conrnkap}
 \min_{X\in \cR(p,q,k)} ||A -BXC||_F,
 \end{equation}
 having the minimal $||X||_F$.  This solution is unique if and only if either
 $$
 k\ge \rank P_{B,L}AP_{C,R}\quad \mbox{or} \quad 1\le k < \rank P_{B,L}AP_{C,R}
 $$ and
 $$
 \sigma_k(P_{B,L}AP_{C,R})>
 \sigma_{k+1}(P_{B,L}AP_{C,R}).
 $$
 \end{theo}

 \textbf{Proof of Theorem \ref{comrnkapt}}
 Recall that the Frobenius norm is invariant under the multiplication
 from the left and the right by the corresponding unitary matrices.
 Hence $||A-BXC||_F=||\widetilde {A}-\Sigma_B \widetilde {X} \Sigma_C||$, where
 $\widetilde {A}:=U_B^* A V_C$ and  $\widetilde {X}:=V_B^* X U_C$.
 Clearly, $X$ and $\widetilde X$ have the same rank and the same Frobenius norm.
 Thus, it is enough to consider the minimal problem
 $\displaystyle{\min_{\widetilde {X}\in \cR(p,q,k)}
 ||\widetilde {A} -\Sigma_B\widetilde {X} \Sigma_C||_F}$.

 Let $s=\rank B$ and  $t=\rank C$.  Clearly if $B$ or $C$ is a zero matrix, then $X=\0$
 is the solution to the minimal problem (\ref{conrnkap}).  In this case either $P_{B,L}$ or
 $P_{C,R}$ are zero matrices, and the theorem holds trivially in this case.

 Let us consider the case
 $1\le s, 1\le t$.  Define $B_1:=\diag(\sigma_1(B),\ldots,\sigma_s(B))\in \C^{s\times s},
 C_1:=\diag(\sigma_1(C),\ldots,\sigma_t(C))\in \C^{t\times t}$.
 Partition $\widetilde {A}$ and
 $\widetilde {X}$ into four block matrices $ A_{ij}$ and $ X_{ij}$ with $i,j=1,2$
 so that $\widetilde {A}=[A_{ij}]_{i,j=1}^2$ and $\widetilde {X}=[X_{ij}]_{i,j=1}^{2}$,
 where $A_{11},X_{11}\in \C^{s\times t}$.
 (For certain values of $s$ and $t$, we may have to
 partition $\widetilde A$ or $\widetilde X$ to less than four block matrices.)
 Next, observe  that $Z:=\Sigma_B \tilde X \Sigma_C=[Z_{ij}]_{i,j=1}^2$, where
 $Z_{11}=B_1 X_{11}C_1$ and all other blocks $Z_{ij}$ are zero matrices.
Since $B_1$ and $C_1$ are invertible we deduce
 $$\rank Z=\rank Z_{11}=\rank X_{11}\le \rank \widetilde {X}\le k.$$
The approximation property of $(A_{11})_k$ yields the inequality
 $||A_{11}-Z_{11}||_F\ge ||A_{11}-(A_{11})_k||_F$ for any $Z_{11}$
 of rank $k$ at most.
 Hence for any $Z$ of the above form,
 $$||\widetilde {A} - Z||_F^2= ||A_{11}-Z_{11}||_F^2 + \sum_{2<i+j\le 4} ||A_{ij}||_F^2 \ge
 ||A_{11}-(A_{11})_k||_F^2 + \sum_{2<i+j\le 4} ||A_{ij}||_F^2.$$
 Thus $\widehat {X}=[X_{ij}]_{i,j=1}^2$, where $X_{11}=B_1^{-1}(A_{11})_k C_1^{-1}$ and $X_{ij}=\0$
 for all $(i,j) \ne (1,1)$ is a solution to the problem
 ${\displaystyle \min_{\widetilde {X}\in \cR(p,q,k)}
 ||\widetilde {A} -\Sigma_B\widetilde {X} \Sigma_C||_F}$
 with the minimal Frobenius form.
 This solution is unique if and only if the solution
 $Z_{11}=(A_{11})_k $ is the unique solution to the problem
 ${\displaystyle \min_{Z_{11}\in \cR(s,t,k)} ||A_{11}-Z_{11}||_F}$.
 This happens if either $k\ge \rank A_{11}$ or $1\le k < \rank A_{11}$ and
 $\sigma_k(A_{11}) > \sigma_{k+1}(A_{11})$.
 A straightforward calculation shows
 that $\widehat {X}= \Sigma_B^{\dagger}(P_{\Sigma_B,L} \widetilde {A}
 P_{\Sigma_C,R})_k \Sigma_C^{\dagger}$.
 Thus, a solution of (\ref{conrnkap}) with the minimal Frobenius norm is given by
 \begin{eqnarray*}
 X&=&B^{\dagger}U_B(P_{\Sigma_B,L} U_B^* A V_C P_{\Sigma_C,R})_k  V_C^* C^{\dagger}\\
 &=&B^{\dagger}U_B(U_B^*P_{B,L} A P_{C,R}V_C)_k  V_C^* C^{\dagger}\\
 &=&B^{\dagger}(P_{B,L} A P_{C,R})_k C^{\dagger}.
\end{eqnarray*}

 This solution is unique if and only if
 either $k\ge \rank P_{B,L}AP_{C,R}$
 or $1\le k < \rank P_{B,L}AP_{C,R}$ and
 $\sigma_k(P_{B,L}AP_{C,R})>
 \sigma_{k+1}(P_{B,L}AP_{C,R})$.
 \qed

 A special case of the minimal problem (\ref{conrnkap}), where $X$
 is a rank one matrix and $C$ the identity matrix, was considered by
 Michael Elad \cite{ela} in the context of image processing.

 \section{Examples}

 First observe that the classical approximation problem given by
 (\ref{svdap}) is equivalent to the case $m=p, n=q, B=I_m, C=I_n$.
 (Here, $I_m$ is the
 $m\times m$ identity matrix.)
 Clearly
 $P_{I_m,L}=I_m,\; P_{I_n,R}=I_n,\; I_m^{\dagger}=I_m,\; I_n^{\dagger}=I_n$.
 In this case we obtain the classical solution $B^{\dagger} (P_{B,L}A P_{C,R})_k
 C^{\dagger}=A_k$.

 Second, if $p=m,$ $q=n$ and $B$, $C$ are non-singular, then $\rank (BXC) = \rank X.$
 In this case,
 $P_{B,L} = I_m$ and $P_{C,R} = I_n,$ and the solution to (\ref{conrnkap}) is given by
 $X = B^{-1}A_k C^{-1}.$

 Next, a particular case of the problem (\ref{conrnkap}) occurs in study of a random vector
 estimation (see, for example, \cite{yam2,hua1}) as follows.
 Let $(\Omega, \Sigma, \mu)$ be a probability
 space, where $\Omega$ is the set of outcomes, $\Sigma$ a $\sigma$--field
 of measurable subsets $\Delta\subset \Omega$ and $\mu:\Sigma \mapsto [0,1]$
 an associated probability measure on $\Sigma$ with $\mu(\Omega) = 1$.  Suppose that
 $\x \in L^{2}(\Omega,{\mathbb R}^{m})$ and $\y \in L^{2}(\Omega,{\mathbb R}^{n})$
 are random vectors such that $\x = (x_1,\ldots,$ $ x_m)^T$ and $\y = (y_1,\ldots,
 y_n)^T$ with
 $x_i, y_j\in L^{2}(\Omega,{\mathbb R})$ for  $i=1,\ldots,m$ and $j=1,$ $\ldots,n$,
 respectively.
 Let $E_{xy}=[e_{ij,xy}]\in\R^{m\times n}, E_{yy}=[e_{jk,yy}]
 \in\R^{n\times n}$ be correlation matrices with entries

 \begin{eqnarray*}
 e_{ij,xy}= \int_\Omega x_i(\omega)y_j(\omega)d\mu(\omega),\quad
 e_{jk,yy}= \int_\Omega y_j(\omega)y_k(\omega)d\mu(\omega),\\
 i = 1,\ldots,m, \quad j,k = 1,\ldots,n, \quad \omega\in \Omega.
 \end{eqnarray*}
 The problems
 considered in \cite{yam2,hua1} are reduced to finding a solution to
 the problem (\ref{conrnkap}) with $A = E_{xy} E_{yy}^{1/2 \dag},$
 $B=I_n$ and $C = E_{yy}^{1/2}$
 % ${\displaystyle \min_{X\in \in \cR(m,n,k)} ||E_{xy} E_{yy}^{1/2 \dag}
 %- XE_{yy}^{1/2}||_F,}$
 where we write $E_{yy}^{1/2 \dag}= (E_{yy}^{1/2})^\dag.$ Let the SVD
 of $E_{yy}^{1/2}$ be given by $E_{yy}^{1/2}=V_n\Sigma V_n^*$ and let
 $\rank E_{yy}^{1/2}=r$. Here, $V_n=[\v_1,\ldots, \v_n]$ with $\v_i$ the
 $i$-th column of $V_n.$ By Theorem \ref{comrnkapt}, the solution to
 this particular case of the  problem (\ref{conrnkap}) having the
 minimal Frobenius norm is given by $X = (E_{xy} E_{yy}^{1/2 \dag}V_r
 V_r^*)_k E_{yy}^{1/2 \dag},$ where $E_{yy}^{1/2 \dag}V_r V_r^* =
 E_{yy}^{1/2 \dag}.$ Therefore, $X = (E_{xy} E_{yy}^{1/2 \dag})_k
 E_{yy}^{1/2 \dag}.$ The conditions for the uniqueness follow
 directly from Theorem \ref{comrnkapt}.

\end{document}